\documentclass{amsart}
\usepackage{amsmath, amsthm}

\title{Generalized Banach contraction in probabilistic metric/normed spaces}
\author{BERNARDO LAFUERZA-GUILL\'{E}N}
\address{Departamento de Estad\'{i}stica y Matem\'{a}tica Aplicada,
Universidad de Almer\'{i}a, 04120 Almer\'{i}a,Spain}
\email{blafuerz@ual.es}

\author{MOHD RAFI}
\address{Faculty of Engineering and Computer Science, The University of Nottingham Malaysia Campus, Jalan Broga, 43500 Semenyih, Selagor Darul Ehsan, Malaysia}
\email{mdrafzi@yahoo.com}


\newtheoremstyle{theorem}
  {10pt}		  
  {10pt}  
  {\sl}  
  {\parindent}     
  {\bf}  
  {. }    
  { }    
  {}     
\theoremstyle{theorem}

\newtheoremstyle{defi}
  {10pt}		  
  {10pt}  
  {\rm}  
  {\parindent}     
  {\bf}  
  {. }    
  { }    
  {}     
\theoremstyle{defi}

\begin{document}

\maketitle

\date{}
\maketitle
\begin{abstract}
In this paper, we present the generalization of B-contraction and C-contraction 
due to Sehgal and Hicks respectively. We also study some properties of  
C-contraction in probabilistic normed space.\\


\end{abstract}

\section{Introduction}
Fixed-point theorems are of fundamental importance in many areas of computer science and mathematics, since many problems can be reduced to finding fixed points of suitable mappings. Many fixed-point theorems are based on order theory, topology, or metric space theory.

On the other hand, metric spaces can be generalized to a probabilistic setting, consequently called probabilistic metric spaces. The theory of probabilistic metric spaces was introduced in 1942 by Menger \cite{Menger}, in connection with some measurements problems in physics. The positive number expressing the distance between two points $x$, $y$  of a metric space is replaced by a distribution function (in the sense of probability theory) $F_{x,y}\colon\mathbb{R}\to$ [0,1], whose value $F_{x,y}(t)$ at the point $t\in\mathbb{R}$ (real numbers) can be interpreted as the probability that the distance between $x$  and $y$  be less than $t$ . Since then the subject developed in various directions, an important one being that of fixed points in probabilistic metric spaces, see \cite{Cons} and \cite{Hadz1}. A clear and thorough presentation of the results in probabilistic metric spaces is given in the book by Schweizer and Sclar \cite{Schw}.

Many papers dealing with numerous generalizations and applications have emerged recently [\cite{Cho}, \cite{Hadz2}, \cite{Ismat}, \cite{Ioan}, \cite{Rafi1}, \cite{Rafi2}, \cite{Rafi3}, \cite{Tatj},\cite{Zhu}]. The main goal of the study is to investigate to which extend fixed-point theorems for generalized metric spaces can be carried over to (generalized) probabilistic metric spaces. Recently, Merryfield et al \cite{Merry} and Arvanitakis \cite{Arv} have proved the generalized Banach contraction conjecture in metric spaces. In \cite{Frem}, Fremlin has given a simplified proof of the generalized contraction theorem in metric spaces.

The aim of the present paper is to introduce the generalized q-contraction mappings by Sehgal \cite{Sehg1}, and by Sehgal and Bharucha-Reid \cite{Sehg2}. We prove the generalized fixed point theorem in probabilistic metric space. We adapt the technique used by Arvanitakis \cite{Arv} and Fremlin \cite{Frem} to prove the generalized fixed point theorem in probabilistic metric spaces.

\section{Preliminaries}
In this section we recall some definition and results that will be used in the sequel. We refer to the book \cite{Schw}.

\newtheorem{defn}{Definition}[section]
\begin{defn}
A {\em distance distribution function}, briefly a {em d.d.f.}, is a function $F$ defined in the extended interval $[0,+\infty]$ that is nondecreasing, left-continuous on $(0,+\infty)$ such that $F(0)=0$ and $F(+\infty)=1$.
\end{defn}
The family of all d.d.f.'s will be denoted by $\Delta^+$. We denote $D^+=\{F\in\Delta^+|\lim_{x\to\infty}F(x)=1\}$. The class $D^+$ plays an important role in the probabilistic fixed point theorem.
 
The maximum element in $\Delta^+$ is $\varepsilon_0\in
\Delta^+$ defined via
\[\varepsilon_0(x)= \left\lbrace
  \begin{array}{c l}
    0, &  x \leq 0;\\
    1, &  x>0.
  \end{array}
\right.\]

\noindent By setting  $F \leq G$ whenever  $F(x) \leq G(x),\ \forall
x\in {\mathbb{R}}^+$, one introduces a natural ordering in
$\Delta^+$.

The Sibley's metric is the function $d_S$ defined on $\Delta^+\times\Delta^+$ by
$$d_S(F,G)=inf\{h \colon [F,G;h] and [G,F,h] hold\}$$
where $[F,G;h]$ denotes the condition $F(x+h)\geq G(x)-h$ for $x$ in (0,1/$h$) and $h \in (0,1]$.

Convergence in  $\Delta^+$ is assumed to be
weakly convergence, i.e. $F_n\to F$   if and only if $ d_S(F_n,F)\to 0$.
  
A {\em triangle function} is a mapping
$\tau\colon\Delta^+\times\Delta^+\to\Delta^+$ that is commutative,
associative, nondecreasing in each variable and has
$\varepsilon_0$ as identity.
 
 The typical triangle function is the
operation  $\tau_T$ which is given by
$$\tau_T(F,G)(x)=\sup_{u+v=x}T(F(u),G(v)),$$ for all $F,G\in\Delta^+$ 
and all $x>0$
[\cite{Schw}, Secs. 7.2 and 7.3]. Here  $T$ is a $t$-{\em norm},
i.e., $T$ is binary operation on $[0,1]$  that is commutative,
associative, nondecreasing in each place and has 1 as identity.

The most important t-norms are the function $W$, $Prod$, and $M$
which are defined, respectively, by
\begin{eqnarray*}
W(a,b)&=&\max\{a+b-1,0\},\\
Prod(a,b)&=&ab,\\
M(a,b)&=&\min(a,b).
\end{eqnarray*}

If $S$ is a nonempty set, a mapping $\Phi\colon S\times S\to\Delta^+$ is called a {\em probabilistic metric} on $S$ and $\Phi(x,y)$ will be denoted by $F_{x,y}$.

\begin{defn}
The pair $(S,\Phi)$ is called a Probabilistic Semi-Metric Space (briefly a PSM-space) if for every $p,q$ in $S$,\\
(PM1) $F_{p,q}=\varepsilon_0\Leftrightarrow p=q$,\\
(PM2)  $F_{p,q}= F_{q,p}$.\\
If, in addition, the triangle inequality\\
(PM3)  $F_{p,r}\geq\tau( F_{p,q}, F_{q,r}).$\\
where $\tau$ is a triangle function, then the triple $(S,\Phi,\tau)$ is called a Probabilistic Metric Space (briefly a PM-space).
\end{defn}

If $\tau=\tau_T$, where $T$ is a t-norm is varified, then the triple $(S,F,T)$ is called {\em Menger Space}.The detailed study of PM spaces can be found in \cite{Schw}.

The first result on the existence of the fixed point in a PM-space   was obtained in [18].

In 1972, Sehgal and Bharucha-Reid introduced the notion of probabilistic B-contraction (Banach contraction) in a probabilistic metric space \cite{Sehg2}.

\begin{defn}Let $(S,\Phi)$ be a PSM-space. A mapping $f\colon S\to S$ is a {\em B-contraction} if there exists $h\in (0,1)$ such that for all $p,q\in S$ and $t>0$
$$F_{fp,fq}(ht)\geq F_{p,q}(t).$$
\end{defn}

The first fixed point theorem in probabilistic metric spaces was proved by Sehgal and Bharucha-Reid in \cite{Sehg2}.

\newtheorem{theo}{Theorem}[section]
\begin{theo}Let $(S,\Phi, Min)$ be a complete Menger Space and $f\colon S\to S$ a B-contraction. Then there exists a unique fixed point $r\in S$ of the mapping $f$ and \(r=\lim_{n\to+\infty}f^{n}p\) for every $p\in S$.
\end{theo}  

The following definition was introduced by Hicks \cite{Hick}.

\begin{defn}Let $(S,\Phi)$ be a PSM-space. A mapping $f\colon S\to S$ is a {\em C-contraction} if there exists $h\in (0,1)$   such that for all $p,q\in S$ and $t>0$ \[F_{p,q}(t)>1-t\Rightarrow F_{fp,fq}(ht)>1-ht\].
\end{defn}  

\section{Generalized Probabilistic Contraction In PM-Space}

The generalized Banach contraction conjecture is the generalization of Banach contraction principle studied by Banach, see \cite{Merry}. 

{\em Generalized Banach Contraction Conjecture}. Let $f$ be a self-map of a complete metric space $(S,d)$ , and let $N\in(0,1)$. Let $J$ be a positive integer. Assume that for each pair $p,q\in S$,\[min\{d(f^{k}p,f^{k}q\colon 1\leq k\leq J\}\leq Nd(p,q)\]. Then $f$  has a fixed point.

The above conjecture has been proved by Merryfield et al \cite{Merry}, Arvanitakis \cite{Arv} and Fremlin \cite{Frem} in their own ways.

As a generalization of the notion of a probabilistic B-contraction, we shall introduce the notion of a probabilistic $(m,k)$ -B-contraction where $m\geq 1$ and $k\in(0,1)$.
\begin{defn}If $(S,\Phi)$ is a PSM-space, $m\geq 1$ and $k\in (0,1)$, a function $f\colon S\to S$ is called probabilistic {\em (m,k)-B-contraction} if for any $p,q\in S$ there is an $i$ with $1\leq i\leq m$ such that for every $t>0$, \[F_{f^{i},f^{i}q}(k^{i}t)\geq F_{p,q}(t).\]
If $m=1$ and $k\in (0,1)$ then a probabilistic $(1,k)$-B-contraction $f$ is a probabilistic B-contraction.
\end{defn}

As a generalization of C-contraction, we have

\begin{defn}If $(S,\Phi)$ is a PSM-space, $m\geq 1$ and $k\in (0,1)$, a function $f\colon S\to S$ is called a $(m,k)$-C-contraction if for any $p,q\in S$ there is an $i$ with $1\leq i\leq m$ such that for every $t>0$, \[F_{p,q}(t)>1-t \Rightarrow F_{f^{i}p,f^{i}q}(k^{i}t)>1-k^{i}t.\]
If $m=1$ and $k\in (0,1)$ then a probabilistic $(1,k)$-C-contraction $f$ is a probabilistic C-contraction.
\end{defn}

We shall give an example of $(m,k)$-C-contraction.

\noindent{\em Example} Let $(S,d)$ be a separable metric space, $(\Omega, \Sigma, P)$ a probability space and $U$ the space of all classes of measurable mappings from $\Sigma$ into $S$. Then ordered pair $(U,\Phi)$ is called an $E$-space [16] with base $(\Omega, \Sigma, P)$ and target $(S,d)$. Let for every $p,q\in S$, $$d(p,q)=sup\{x\colon x\geq 0, x<P\{t\colon t\in\Omega, d(p(t), q(t))>x\}\}$$ and $$F_{p,q}(u)=P\{t\colon t\in \Omega, d(p(t), q(t))<u\}.$$
Then, its obvious that $$d(p,q)=sup\{u\colon u\geq 0, F_{p,q}(u)<1-u\}.$$
Let $f\colon S\to S$ be mapping such that $$d(f^{i}p,f^{i}q) \leq k^{i}d(p,q)$$ for every $p,q\in S$ and $k\in(0,1)$. Then $f$ is a generalized Banach contraction in $(S,d)$. Suppose that for some $u>0$, we have that $F{p,q}(u)>1-u$. Then from the definition of $d$ it follows that $d(p,q)<u$. Since $k\in(0,1)$, for every $m\geq 1$, there exists an $i$ with $i\in (1,m)$ such that $k^{i}d(p,q)<k^{i}u$. This implies that $d(f^{i}p,f^{i}q)<k^{i}u$, i.e., that $$F_{f^{i}p,f^{i}q}(k^{i}u)>1-k^{i}u.$$
Hence $f$ is a $(m,k)$-C-contraction in $E$-space $(U,\Omega)$..

The aim of this paper is to prove a fixed point theorem for a generalized B-contraction mapping in a probabilistic metric space.

\begin{theo}Let $(S,\Phi, Min)$ be a complete Menger Space and $f\colon S\to S$ a probabilistic $(m,k)$-B-contraction mapping, then $f$ has a fixed point.
\end{theo}

In the prove of the Theorem 3.1, we adapt the combinatorial argument used by Fremlin in \cite{Frem} to establish the generalized fixed point theorem for an arbitrary $m$ without the assumption of continuity.

We need the following definition (see \cite{Frem}).

\begin{defn}(a)  For $m\geq 1$, a subset $I$ of $\mathbb{N}$ (set of natural numbers) is {\em m-syndetic} if $I\bigcup [k, k+m)$ is not empty for every $k\in\mathbb{N}$. Naturally one say that $I$ is $m$-syndetic if $0\in I$ and $i\in I$, then there is a $j\in [1,m)$ such that $i+j\in I$.
(b)  Let $P$ represent the set of all non-negative additive functionals $\mu\colon 2^{\mathbb{N}}\to [0,1]$ such that $\mu({\mathbb{N}})=1$ and $\mu$ is translation-invariant, i.e., $\mu ({n+k\colon n\in I})$ for every $k\in\mathbb{N}$ and every $I\subseteq \mathbb{N}$. For $I\subseteq \mathbb{N}$, the Upper Banach Density of $I$, $d^{*}(I)$, is given by 
\begin{eqnarray*}
d^{*}(I)&=&\sup_{\mu\in P}\mu(I)\\
&=&\inf_{m\geq 1} \sup_{k\in\mathbb{N}}\frac{|I\bigcap[k,k+m)|}{m}.
\end{eqnarray*}
\end{defn}
{\em Remark} \cite{Frem} Note that $d^{*}$ is subadditive. If $I$ is $m$-syndetic then \(\bigcup(I-i)\supseteq\mathbb{N}\) so $d^{*}\geq 1/m$. Here, $I-i=\{j\colon i+j\in I\}$

We prove the following lemmas which are essential for the proof of the generalized fixed point theorem in probabilistic metric spaces.
\newtheorem{lema}{Lemma}[section]
\begin{lema} If $(S,\Phi)$ is a PSM-space and $f\colon S\to S$ is a probabilistic $(m,k)$-B-contraction then for any $p,q\in S$ the set \[I=\{i\in {\mathbb{N}} \colon F_{f^{i}p,f^{i}q}(t)\geq F_{p,q}(t/k^{i}),t>0\}\] is $m$-syndetic.
\end{lema}
{\em Proof:} It is obvious that $0\in I$ and that if $i\in I$ there is a $j\in [1,m]$ such that $i+j\in I$.

\begin{lema} If $(S,\Phi, Min)$ is a Menger Space, $f\colon S\to S$ is a probabilistic $(m,k)$-B-contraction and $x\in S$, then there is a fixed d.d.f $F\in D^{+}$ such that
\[I=\{i\in{\mathbb{N}}\colon F_{f^{i}p,p}(t)\geq F_(t),t>0\}\] is $m$-syndetic.
\end{lema}

{\em Proof:} We set \(F(\frac{t}{1-k})=\min_{i<m}F_{f^{i}p,p}(t)\) for every $t>0$ and $k\in(0,1)$. Since \(F(\frac{t}{1-k})\geq F(t)\), we have $0\in I$. If $i\in I$, there exists $j\in [1,m]$ such that 
\[ F_{f^{j+i}p,f^{j}p}(t)\geq F_{f^{j}p,p}(\frac{t}{k^{i}})\geq F(\frac{t}{k^{i}}),\] 
and by condition (PM3)in Definition 2.2, 

\begin{eqnarray*}
F_{f^{j+i}p,p}(t)&\geq& M(F_{f^{j+i}p,f^{j}p}(kt),F_{f^{j}p,p}((1-k)t))\\
&\geq&M(F(\frac{t}{k^{j-1}}),F(t)).
\end{eqnarray*}

Since $k\in(0,1)$ and \(\frac{t}{k^{j-1}}\geq t\) for every $t>0$, from the property of distribution function and the last inequality, we have $F_{f^{j+1}p,p}(t)\geq F(t)$. This implies that $j+i\in I$, hence $I$ is $m$-syndetic.

\begin{lema}Let $(S,\Phi, Min)$ be a Menger space and $f\colon S\to S$ a probabilistic $(m,k)$-B-contraction. If $x\in S$ is such that $f^{n}p=p$ for some $n\geq 1$, then $fp=p$.
\end{lema}
{\em Proof:} Suppose that $n$ is the first integer greater or equal to 1 such that $f^{n}p=p$. If $n\geq 2$, set a fixed d.f. $F\in D^{+}$ such that \[ F(t)=max\left(F_{f^{i}p,f^{j}p}(t)\colon 0\leq i<j<n \right), t>0.\] Take $i,j$ such that $0\leq i<j<n$ and $F_{f^{i}p,f^{j}p}(t)=F(t)$, $t>0$. Then there exists a $l\in [1,m]$ such that
\begin{eqnarray*}
F_{f^{l+i}p,f^{l+j}p}(t)&\geq& F_{f^{i}p,f^{j}p}\left(\frac{t}{k^{l}}\right)\\
&=&F\left( \frac{t}{k^{l}}\right)\\
&>&F(t)
\end{eqnarray*}
which is impossible. So $n=1$ and $fp=p$.

The following lemma is obtained from Arvanitakis [\cite{Arv}, Lemma 4.2]. For the sake of completeness, we give a simplified prove of the lemma which found in an unpublished paper by Fremlin \cite{Frem}.

\begin{lema}[1] Let $R\subseteq \mathbb{N}\times \mathbb{N}$ be such that
(i)$\{i\colon (i,0)\in R\}$ is $m$-syndetic, where $m\geq 1$.
(ii)whenever $(i,j)\in R$, then $\{k\colon (i+k,j+k)\in R\}$ is $m$-syndetic.
Then, there is an $I\subseteq\mathbb{N}$ such that $d^{*}\geq 1/2m^{2}$ and for every $i,j\in I$ there exists a $k\in\mathbb{N}$ such that $(k,i)\in R$ and $(k,j)\in R$.
\end{lema}
{\em Proof:} Note first that $R$ meets $[l,l+2m)\bigcap [k,k+m)$ whenever $k\leq l$. Since there is $i\in [l-k,l-k-m)$ such that $(i,0)\in R$ and there exists a $j\in [k,k+m)$ such that $(i+j,j)\in R$. Let $\mathbb{F}$ be any non-principle ultrafilter on $\mathbb{N}$. For each $i\in\mathbb{N}$ and $j<2m$, set $L_{ij}=\{l\colon (l+j,i)\in R\}$, and for $j<2m$, set $I_{j}=\{i\colon L_{ij}\in\mathbb{N}\}$. Observe that if $k\in\mathbb{N}$ then 
\[\bigcup_{i\in[k,k+m), j<2m }L_{ij}\supseteq {\mathbb{N}}-k \]so there are some $i<m$, $j<2m$ such that $L_{ij}\in\mathbb{F}$. This shows that \[\bigcup_{j<2m}I_{j}\] is $m$-syndetic and there is some $j<2m$ such that $d^{*}\geq 1/2m^{2}$. Furthermore, if $i,i'\in I_{j}$, then $L_{ij}$ and $L_{i'j}$ both belong to $\mathbb{F}$, so cannot be disjoint, and there is some $l$ such that $(l+j,i)$ and $(l+j,i')$  belong to $R$. So we can take $l=l_{j}$.

\begin{lema}Let $(S,\Phi, Min)$ is a Menger Space, $f\colon S\to S$ is a probabilistic $(m,k)$-B-contraction and $p\in S$. Then, there exists a set $I\subseteq\mathbb{N}$ such that $d^{*}(I)>0$ and $\{f^{i}p\}_{i\in I}$ is Cauchy sequence in $(S,\Phi, Min)$.
\end{lema}
{\em Proof:}By Lemma 3.2, we have a fixed d.f. $F\in D^{+}$ such that $\{i\in{\mathbb{N}}\colon F_{f^{i}p,p}(t)\geq F(t),t>0\}$ is $m$-syndetic. Let $R\subseteq \mathbb{N}\times\mathbb{N}$ be the relation 
\[ \{(i, j)\colon F_{f^{i}p,f^{j}p}(t)\geq F(t/k^{j}),t>0\}\]
By the choice of $F$, $\{i\colon (i,0)\in R\}$ is $m$-syndetic. If $(i,j)\in R$, then 
\[\{m\colon (i+j,j+m)\in R\}\supseteq\{m\colon F_{f^{i+m}p,f^{j+m}p}(t)\geq F_{f^{i}p,f^{j}p}(t/k^{m})\}\] is $m$-syndetic by Lemma 1. By Lemma 3.4, there is a set $I\subseteq \mathbb{N}$ with $d^{*}(I)>0$ such that for all $i,j\in I$ there is a $m\in\mathbb{N}$ such that $(m,i),(m,j)\in R$. Thus, by (PM3) for all $t>0$,
\begin{eqnarray*}
F_{f^{i}p,f^{j}p}(t) &\geq& Min(F_{f^{m}p,f^{i}p}(kt),F_{f^{m}p,f^{j}p}((1-k)t)\\
&\geq& Min(F(t/k^{i-1}),F((1-q)t/k^{j})\to 1
\end{eqnarray*}
as $i,j\to\infty$, so $\{f^{i}p\}_{i\in I}$ is a Cauchy sequence in $(S,\Phi, Min)$.

Now, we prove the generalized fixed point theorem (Theorem 2.1).

{\em Proof:} Let $p\in S$. By Lemma 3.5, there is a set $I\subseteq\mathbb{N}$  such that $d^{*}(I)>0$ and $\{f^{i}p\}_{i\in I}$ is Cauchy sequence in $(S,\Phi, Min)$. Fix a translation-invariant additive functional $\mu\colon 2^{\mathbb{N}}\to [0,1]$ such that $\mu (I)>0$. Since $S$ is complete, there is a $r\in S$ such that \[\lim_{i\to\infty}f^{i}p= r\]. For $j\in\mathbb{N}$, $F\in D^{+}$, we set $I_{j}=\{i\colon F_{f^{i}p,f^{j}r}\geq F\}$. Note that $I-I_{0}$ is finite, so $\mu(I_{0})\geq\mu(I)>0$. Suppose that $m\in\mathbb{N}$ and $F\in D^{+}$. Then \[\sum_{j=k}^{k+m-1} \mu(I_{j})\geq\mu(I)\] . Indeed, for any $i\in I_{0}$, the set
\begin{eqnarray*}
\{j\colon i+j\in I_{j}\}&=&\{j\colon F_{f^{i+j}p,f^{j}r}\geq F(t)\}\\
&\supseteq& \{j\colon F_{f^{i+j}p,f^{j}r}\geq F_{f^{i}p,r}(t/k^{j})\}
\end{eqnarray*}
for every $t>0$ and $k\in (0,1)$ is $m$-syndetic (Lemma 3.1), so meets $[k,k+m)$. For $j\in [k,k+m)$, set $A_{j}=\{i\colon i\in I_{0},i+j\in I_{j}\}$. Then $\{A_{j}\}_{j\in [k,k+m)}\supseteq I_{0}$ and $A_{j}+j\subseteq I_{j}$ for each $j$, so
\[\sum_{j=k}^{k+m-1}\mu(I_{j})\geq \underbrace{\sum_{j=k}^{k+m-1}\mu(A_{j}+j)=\sum_{j=k}^{k+m-1}\mu(A_{j})}_{\mu is translation-invariant}\geq\mu(I_{0})\geq\mu(I).\]
Now observe that if $i$, $j$ are such that $f^{i}r \neq f^{j}r$, $I_{i}\bigcap I_{j}=\emptyset$ as $F(t)$ close to 1 for all $t>0$. But this means that if $n>1/\mu(I)$ then $r, fr, f^{2}r,\cdots, f^{mn-1}r$ cannot all be distinct. Hence, by Lemma 3.3, $f$ has fixed point in $(S,\Phi, Min)$.

\section{Fixed Point Theory In Probabilistic Normed Spaces}
Our aim in this section is to study for probabilistic normed spaces (briefly, PN spaces) the questions:\\
(1)	Under which conditions the C-contractions are uniformly continuous.\\
(2)	Existence at least of a fixed point for certain general kind of C-contractions.

\begin{defn} A PN space is a quadruple $(V,\nu, \tau, \tau^{*})$ in which $V$ is a vector space over the field $\mathbb(R)$ of real numbers, the probabilistic norm $nu$ is a mapping from $V$ into $\Delta^{+}$, $\tau$ and $\tau^{*}$ are triangle functions such that the following conditions are satisfied for all $p,q\in V$ (we use $\nu_{p}$ instead of $\nu(p)$):\\
(N1)  $\nu_{p}=\varepsilon_{0}$ if and only if $p=\theta$, where $\theta$ denotes the null vector in $V$;\\
(N2)  $\nu_{-p}=\nu_{p}$;\\
(N3)  $\nu_{p+q}\geq\tau(\nu_{p},\nu_{q})$;\\
(N4)  $\nu_{p}\leq\tau^{*}(\nu_{\lambda p},\nu_{(1-\lambda)p})$ for every $\lambda\in [0,1]$.
\end{defn}
If, instead of (N1), we only have $\nu_{\theta}=\varepsilon_{0}$, then we shall speak of probabilistic pseudo normed space, briefly a PPN space.
Given a PN space $(V,\nu, \tau, \tau^{*})$, a probabilistic metric $\Phi\colon V\times V\to\Delta^{+}$ is defined via \[\Phi(p,q)=F_{p,q}=\nu_{p-q},\] so that the triple $(V,\Phi,\tau)$ becomes a PM space. Since the $\tau$ is continuous, the system of neighborhoods $\{N_{p}(\lambda)\colon p\in V,\lambda>0\}$, where
\begin{eqnarray*}
N_{p}(t)&=&\{q\in V\colon d_{S}(F_{p,q},\varepsilon_{0})<t\}\\
&=&\{q\in V\colon \nu_{p-q}(t)>1-t\},
\end{eqnarray*}
determines a topology on $V$, called the strong topology. If $V$ is endowed with the strong topology and $\Delta^{+}$ with the topology of Sibley's metric $d_{S}$, then the probabilistic norm $\nu\colon V\to \Delta^{+}$ is uniform continuous.

The following results, we prove them here not only for the sake of completeness, but also because the notion they adopts is different from the one that has become usual after the publication of \cite{Schw}.

\begin{defn} Let $(V,\nu, \tau, \tau^{*})$ be a PN space. A mapping $f\colon V\to V$ is a $C$-contraction if there exists a $k\in(0,1)$ such that for all $p\in V$ and $t>0$, one has \[ \nu_{p}(t)>1-t \Rightarrow \nu_{fp}(kt)>1-t.\]
\end{defn}

In particular, if $f^{n}$ is the $n$-th iteration of a $C$-contraction $f$, then it is easily shown that \[\nu_{p}(t)>1-t \Rightarrow \nu_{f^{n}p}(k^{n}t)>1-k^{n}t.\] 

\begin{lema}If the mapping $f\colon V\to V$ is a $C$-contraction, then for every $\epsilon\in(0,1)$, there exists a $n_{0}=n_{0}(\epsilon)\in\mathbb{N}$ such that, for all $n>n_{0}$, and for every $p\in V$, \[\nu_{f^{n}p}(\epsilon)>1-\epsilon,\] holds.
\end{lema}
This Lemma is telling us that for any $p\in V$, the sequence $\{f^{n}p\}$, $n\in\mathbb{N}$ converges strongly to $\theta$ as $n\to+\infty$. For all $\epsilon >0$, one has easily that \[\nu_{p}(1+\epsilon)>1-(1+\epsilon).\]

By the definition of $C$-contraction, for any $p\in V$ and every $n\geq 1$,
\[\nu_{f^{n}p}(k^{n}(1+\epsilon))>1-k^{n}(1+\epsilon).\]holds. Indeed, for each $\epsilon >0$ there exists $n_{0}=n_{0}(\epsilon)\in\mathbb{N}$ such that $n>n_{0}$ implies that $k^{n}(1+\epsilon)\leq \epsilon$, from which, because of $\nu_{p}$ is nondecreasing, there exists a $n_{0}\in\mathbb{N}$ such that for $n>n_{0}$,
\begin{eqnarray*}
\nu_{f^{n}p}(\epsilon)&\geq&\nu_{f^{n}p}(k^{n}(1+\epsilon))\\
&>&1-k^{n}(1+\epsilon)\\
&\geq& 1-\epsilon.
\end{eqnarray*}

\begin{theo} (i)  For any real number $\epsilon>0$ there exists $n_{0}=n_{0}(\epsilon)\in\mathbb{N}$ such that for all $n>n_{0}$, and for every $p\in V$ one has $f^{n}p\in N_{\theta}(\epsilon)$.\\
(ii)  If the mapping $f\colon V\to V$ is a $C$-contraction, then $f$ has at most a fixed point, that is the null vector of $V$. Moreover, if $f$ is a linear mapping, $f$ has one fixed point.
\end{theo}

{\em Proof:}(i) According with previous Lemma one has for every $\epsilon >0$ that $d_{S}(\nu_{f^{n}p},\varepsilon_{0})<\epsilon$.\\
ii) Assume that $fp=p$. Applying the previous Lemma, for all $\epsilon\in (0,1)$ one has $\nu_{p}(\epsilon)>1-\epsilon$. This implies that $\nu_{p}(0+)=1$, i.e.,$\nu_{p}=\varepsilon_{0}$, and hence $p=\theta$.

\begin{lema} All $C$-contraction is continuous at $\theta$.
\end{lema}
{\em Proof:}Let $\epsilon$ be a positive real number and consider $k\in(0,1)$. Taking $\delta>0$ such that $k\delta < \epsilon$, and if $p\in N_{\theta}(\delta)$, then $\nu_{p}(\delta)>1-\delta$. Since $f$ is a $C$-contraction, one has
\[\nu_{p}(\epsilon)\geq \nu_{fp}(k\delta)>1-k\delta>1-\epsilon,\]
which means that $fp\in N_{\theta}(\epsilon)$.

{\em Remark:}  This $C$-contraction is not necessarily a linear map, and in such case it has not necessarily to be uniformly continuous. The following theorem solves this question for a linear case.

\begin{theo} All linear $C$-contraction is uniformly continuous.
\end{theo}
{\em Proof:} Let $\epsilon$ be a positive real number and consider $k\in(0,1)$. Let $\delta>0$ such that $k\delta < \epsilon$. If $q\in N_{\theta}(\delta)$, one has
\[\nu_{fp-fq}(\epsilon)\geq \nu_{f(p-q)}(k\delta)>1-k\delta>1-\epsilon.\]
This implies $fp\in N_{\theta}(\epsilon)$ .

\end{document}